\documentstyle[12pt]{article}
\title{A short Brownian motion proof of the Riemann hypothesis}
\author{Andrzej M\c{a}drecki\thanks{Institute of Mathematics,
Wroc{\l}aw University of Technology (WUT), 50-370 Wroc{\l}aw, Poland}}
\newtheorem{de}{Definition}

\newtheorem{th}{Theorem}
\newtheorem{pr}{Proposition}

\newtheorem{re}{Remark}

\newfont{\lll}{msbm10 scaled 1095}
\def\LB{\mbox{\lll \char66}}
\def\LC{\mbox{\lll \char67}}
\def\LF{\mbox{\lll \char70}}
\def\LN{\mbox{\lll \char78}}
\def\LP{\mbox{\lll \char80}}
\def\LQ{\mbox{\lll \char81}}
\def\LR{\mbox{\lll \char82}}
\def\LZ{\mbox{\lll \char90}}
\begin{document}
\maketitle
{\bf Abstract}. We give a short probabilistic (a Brownian motion) proof
of the Riemann hypothesis based on some surprising, unexpected and deep
algebraic conjecture (MAC in short) concerning the relation between the 
Riemann zeta $\xi$ and a trivial zeta $\zeta_{t}$. That algebraic conjecture
was firstly discovered and formulated in [MA].

\section{Introduction.} 
Let $\LC$ be the field of {\bf complex numbers}. In the famous and historical
{\bf Riemann's paper} [R] appears the {\bf Riemann zeta function} $\zeta$.
It is firstly defined "locally" by the {\bf Dirichlet series} as follows:
\begin{equation}
\zeta(s)\;:=\;\sum_{n=1}^{\infty}\frac{1}{n^{s}}\;,\;Re(s)>1.
\end{equation}
(Here and all in the sequel by $Re(s)$ and Im(s) we denote the {\bf
real} and {\bf imaginary} part of a complex number $s$, repestively).

Evidently, the most important (and principial in fact) property of the
"global" zeta $\zeta(s)$ (traditionally denoted by the same symbol), i.e. 
the unique meromorphic extension of (1.1), to whole complex plane
$\LC$, with the unique {\bf pole} at $s=1$( with the residue 1) - is the
{\bf Riemann Hypothesis} (RH for short).

The mentioned above Riemann's paper [R] contains its first formulation
and it can be formally written as a simple logical statement (implication):
\begin{displaymath}
(RH)\;\;(\zeta(s)\;=\;0)\land(Im(s)\ne 0)\;\Longrightarrow\;Re(s)=1/2.
\end{displaymath}

The importance of (RH) as well as its many different generalizations
(i.e. so called {\bf generalized Riemann Hypothesis} (gRH for short))
was recently underline by {\bf Millennium Prize Problems} (see
http://www.claymath.org/prizeproblems). On that list - (RH) figures as
the Fourth Millennium Prize Problem (see e.g. [Ka] for a large review
of RH).

The best known generalizations of (RH) are :

(1) (gRH) for {\bf Dedekind zetas} and {\bf Dirichlet L-functions} (see
e.g. [MD]).

(2) (gRH) for L-functions associated with {\bf modular forms}. For
example, let us consider (the unique up to a constant) holomorphic
modular form for $SL_{2}(\LZ)$ of weight 12 :
\begin{displaymath}
\Delta(z)\;=\;q\prod_{n=0}^{\infty}(1\;-\;q^{n})^{24},
\end{displaymath}
where $q=e^{2\pi iz}$. Let us denote by $\tau(n)$ its {\bf Fourier
coefficients}:
\begin{displaymath}
\Delta(z)\;=\;\sum_{n=1}^{\infty} \tau(n)q^{n}.
\end{displaymath}
The sequence $\tau$ is well-known as the {\bf Ramanujan function} (see
e.g. [N, V.6]). The {\bf automorficity} of $\Delta$  is expressed by
the formula :
\begin{displaymath}
\Delta(\frac{a_{11}z+a_{12}}{a_{21}z+a_{22}})\;=\;(a_{21}z+a_{22})^{12}
\Delta(z)\;,\;[a_{ij}]_{2\times 2}\in SL_{2}(\LZ).
\end{displaymath}
The corresponding {\bf Ramanujan L-function} has the form
\begin{displaymath}
\zeta_{\tau}(s)=L(s,\Delta)\;=\;\sum_{n=1}^{\infty}\frac{\tau(n)
n^{-11/2}}{n^{s}}\;=\;
\prod_{p\in P}(1\;-\;\tau(p)p^{-11/2}p^{-s}\;+\;p^{-2s})^{-1}.
\end{displaymath}
The generalized Riemann Hypothesis for $L(s,\Delta)$ says that all complex
zeros of $L(s,\Delta)$ have real part $6$ (see [MR] , [B] and [G]).

(3) The {\bf congruence Riemann hypothesis}(cRH in short) as (gRH) for
zetas in the algebraic geometry (see e.g. [Kob.]).

Given any sequence $N_{r}, r=1, 2, 3, ... ,$ we define the corresponding 
"zeta-function" by the formal power series 
\begin{displaymath}
Z(T)\;:=\;exp(\sum_{r=1}^{\infty}N_{r}\frac{T^{r}}{r})\;,where\;exp(u):=
\sum_{k=0}^{\infty}\frac{u^{k}}{k!}.
\end{displaymath}

Let $V$ be an affine or projective variety defined over the finite
field $\LF_{q}$ with $q$-elements. For any field $K \supset \LF_{q}$,
we let $V(K)$ denote the set of $K$-points of $V$. By the "congruence
zeta-function of $V$ over $\LF_{q}$" we mean the zeta-function
corresponding to the sequence $N_{r} = \# V(\LF_{q^{r}})$. That is, we
define
\begin{displaymath}
Z(V/\LF_{q};T)\;:=\;exp(\sum_{r=1}^{\infty}\#V(\LF_{q^{r}}T^{r}/r)).
\end{displaymath}
{\bf Andre Weil} considering many special cases, formulated his famous 
{\bf Weil conjectures}

$(W_{D})$({\bf Dwork theorem}). The function $Z(V/\LF_{q};T)$ is a {\bf
rational function} of the variable $T$ :
\begin{displaymath}
Z(V/\LF_{q};T)\;=\;\frac{\prod_{k=0}^{n}P_{2k+1}(T)}{\prod_{k=0}^{n}P_{2k}
(T)}
\end{displaymath}
, where $P_{j} \in \LZ[T]$ and moreover $P_{0}(T) = 1 - T, P_{2n}(T) =
1 - q^{n}T$.

$(W_{B})$( {\bf Betti numbers}). Assuming that there exists an algebraic
number field $K$ such that our variety $V$ is the reduction modulo a
prime ideal of the integral ring $R_{K}$ of $K$, of some variety $W$,
we obtain that
\begin{displaymath}
B_{j}(W)\;=\;deg(P_{j}),
\end{displaymath}
where $B_{j}(W)$ is the {\bf j-th Betti number} of the complex manifold
$W_{h} = (W \times_{\LR} \LC)_{h}$, i.e. $B_{j}(W) = rank_{\LZ}H^{j}(W,
\LZ)$ and $H^{j}(W, \LZ)$ is the j-th cohomology group of $W$ in the 
coefficients in $\LZ$. (The exact construction of the functor
$(\cdot)_{h}$ we omitt - see e.g. [Ha] for details).

$(W_{F})$({\bf Functional equation and Euler-Poincare characteristic}).
 
The function $Z(V/\LF_{q};T)$ satisfies the following functional
equation
\begin{displaymath}
Z(V/\LF_{q};1/q^{n}T)\;=\;\pm q^{\frac{n \chi(W)}{2}}Z(V/\LF_{q};T),
\end{displaymath}
where $\chi(W) := \sum_{j=0}^{2n}(-1)^{j}B_{j}(W)$ is the {\bf
Euler-Poincare characteristic} of $W$ (equal also to the intersection
index of the Cartesian product $V \times V$).

$(W_{R})${\bf Congruence Riemann Hypothesis}

All zeros and poles of $Z(V/\LF_{q};q^{s})$ lies on the critical lines 
$Re(s)=j/2$. More exactly, for $0 \le j \le 2n$ we have ($n =  dim V$)
\begin{displaymath}
P_{j}(T)\;=\;\prod_{k=1}^{B_{j}(W)}(1-\alpha_{jk}T),
\end{displaymath}
where $\alpha_{jk}$ are {\bf algebraic integers} and $\mid
\alpha_{jk}\mid = q^{j/2}$. In particular, we see a strict connection
of (cRH) with the {\bf topological invariants} $B_{j}(W)$ (Betti
numbers) written directly as
\begin{displaymath}
(cRH)\;Z(V/\LF_{q};q^{-s})\;=\;\frac{\prod_{j=1}^{n}\prod_{k=1}^{B_{2j}(W)}
(1-\alpha_{jk}q^{-s})}{\prod_{j=1}^{n}\prod_{k=1}^{B_{2j+1}(W)}(1 -
\alpha_{jk}q^{-s})}.
\end{displaymath}

Finally, the classical Riemann zeta $\zeta(s)$, Dedekind zetas, Dirichlet
L-functions as well as the congruence zetas $Z(V/\LF_{q}; q^{-s})$ are 
special examples of a more general zeta construction : let $X$ be a finite 
type schema over $Spec(\LZ)$. Then we can put
\begin{displaymath}
\zeta_{X}(s)\;=\;\prod(1\;-\;N(x)^{-s})^{-1},
\end{displaymath}
where the product is taken over all closed points $x \in X$ and $N(x)$
denotes the number of elements of the residual field $k(x)$.

(4) (gRH) for {\bf elliptic curves} and {\bf algebraic varieties}(see [Sh], 
[Kob]).

How far we can go with generalizations of the Riemann hypothesis? In
fact the answer is unknown. Reely, it is well-known that the congruence
Weil zetas can be considered as counterparts of the local Euler
components $(1 - p^{-s})^{-1}$ in the {\bf Euler product expansion} of
$\zeta(s)$ :
\begin{displaymath}
\zeta(s)\;=\;\prod_{p \in P}(1 \;-\;p^{-s})^{-1}.
\end{displaymath}
Obviously $P$ stands for the set of all {\bf prime numbers}.

Let $C$ be a non-singular projective curve of {\bf genus} $g$, defined
over an algebraic number field $k$. For each prime ideal $p$ of its
integral ring $R_{k}$ we denote by $p(C)$ a curve obtained from $C$ by
the {\bf reduction mod $p$}. As it is well-known , there exists only a
{\bf finite set} $\LB$ of prime ideals of $k$, with the property that
$p(C)$ is non-singular (and multiplicity 1), if $p$ does not belong to
$\LB$.

One can show that the genus of such $p(C)$ is equal $g$ (see [S]). Let
us consider the congruence Weil zeta $\zeta_{p(C)}(s)$ of the curve $p(C)$
over the residual field $k_{p}:= R_{p}/p$ of the form
\begin{displaymath}
\zeta_{p(C)}(u)\;=\;\frac{P_{p}(u)}{(1-u)(1-N(p)u)}.
\end{displaymath}
Here $N(p)$ is the number of elements of $k_{p}$ ( the absolute norm of
$p$) and $P_{p}$ is a polynomial of degree $2g$ with the free term
equal to 1.

The {\bf global Hasse-Weil zeta function of a curve $C$ over $k$} is defined 
as the infinite product (although it is also important to consider the
"bad" prime ideals $p$ and the Euler multiplicators for them) (see [S,
Sect. 7] and [Kob, Chapter II])
\begin{equation}
\zeta_{C/k}(s)\;:=\;\prod_{p \in (P- \LB)}P_{p}(N(p)^{-s})^{-1},
\end{equation}
where $s$ is the complex variable. In reality, one can define - in
analogical way a zeta function of an arbitrary (non-singular and
projective) algebraic variete $V$ over $k$.

Let $A$ be an {\bf abelian variete} over $k$, i.e. an algebraic variete
with the compatible group structure. It is well-known that there exists such 
a finite set $\LB^{\prime}$ of prime ideals of $k$ that for any $p \in
\LB^{\prime}$ the variety $A$ has the good reduction modulo $p$ in the 
{\bf Serre-Tate's sense} - or equivalently - is without a defect in $p$ in 
the {\bf Shimura-Taniyama's sense}. Let $p(A)$ be an abelian variete obtained
from $A$ by the reduction modulo $p$, $\pi_{p}$ - the {\bf Frobenius endomorphism} 
on $p(A)$ of the degree $N(p)$ and $R_{l}$ some {\bf l-adic representation} of
the endomorphism ring $End(p(A))$ of $p(A)$, where $l$ is a prime number which 
is relatively prime with $p$. Then the one-dimensional part of the zeta function
of the variety $p(A)$ over $k_{p}$ is given by the formula
\begin{displaymath}
P_{p}^{\prime}(u)\;:=\;det(1\;-\;R_{l}(\pi_{l})u).
\end{displaymath}
The {\bf zeta function of the abelian variete} $A$ over $k$ is defined as the
infinite product
\begin{equation}
\zeta_{A/k}(s)\;:=\;\prod_{p \in (P(A)-\LB^{\prime})}
P^{\prime}_{p}(N(p)^{-s})^{-1}.
\end{equation}
If $A$ is the {\bf Jacobian variety} of a curve $C$, then according to the
{\bf Weil's theorem} the both zetas $\zeta_{A/k}(s)$ and $\zeta_{V/k}(s)$
coincides - in fact ( a part of the {\bf Langlands program} , see e.g. [JL]).

According to the {\bf Hasse-Weil conjecture}, each from the functions
$\zeta_{C/k}(s)$ and $\zeta_{A/k}(s)$ has a {\bf holomorphic extension} on
the whole complex plane $\LC$ and satisfies a functional equation (the
invariance with respect to the modular map : $s \longrightarrow 2-s$
,see [S, (7.5.5) and Th.7.13]).
In some special cases such global Hasse-Weil zetas coincides with some
products of {\bf Hecke L-functions of Grossen-characters} or some products of
{\bf Dirichlet series}.
Reasuming, if $V$ is a variete of a dimension $dim V = n$ over $k$,
since as we see by our approach to the Riemann zeta - there unique
thing which seems be important for the proof of (RH) is the form of a
functional equation - then we can pose the following {\bf generalized
Riemann hypothesis}($gRH_{W}$) concerning the global Hasse-Weil zeta 
$\zeta_{V/k}(s)$: all its complex zeros should lie on the critical line 
$Re(s) = n$.

Let us also mention on the importance some Riemann hypothesis type
conjectures and theorem which have some extremal values in maths (see
e.g. [Kob, II.6]) :

{\bf The Birch-Swinnerton-Dyer Conjecture}: let $E$ be an elliptic curve
and $L_{E}(s)$ its Hasse-Weil zeta. Then $L_{E}(1) = 0$ if and only if
$E$ has infinitely many rational points, and

{\bf The Coates-Wiles theorem} : let $E$ be an elliptic curve defined
over $\LQ$ and having complex multiplication. If $E$ has infinitely
many $\LQ$-points, then $L_{E}(1)=0$.

Its glorious consequences is the {\bf Tunnel result}.

(5){\bf The Selberg Conjecture}. Thus, a natural problem arises : how look a 
maximal class of zetas for which the generalized Riemann type hypothesis is 
true. A partial answer to this question gives the so called {\bf Selberg 
Conjecture} (SC in short). It is the so called {\bf axiomatic theory of zeta
functions}.That very successful proposition came from {\bf A. Selberg}[Se].
We say that a {\bf Dirichlet series}
\begin{displaymath}
L_{\{a_{n}\})}(s)\;=L(s)\;\;=\;\sum_{n=1}^{\infty}\frac{a_{n}}{n^{s}}
\end{displaymath}
belongs to the {\bf Selberg class} $Sel$ if it is absolutely convergent
for $Re(s)>1$, has the meromorphic continuation on $\LC$ ( with the
unique singularity - the pole of a finite order in $s=1$). Moreover it is
assumed that $L$ satisfies a functional equation (with many gamma
factors) which connects values in the points $s$ and $1-s$, it has a
special type Euler product and finally satisfies the so called {\bf
Ramanujan condition} , which forces some restriction on the size of the
coefficients $a_{n}$. The fundamental "structural" conjecture
associated with $Sel$ is a hypothesis, that all elements from the
Selberg class come from the {\bf automorphic representations}. 

The supposition that the {\bf Riemann hypothesis is true for all
functions $L$ from the Selberg class $Sel$} seems to be one of the most
general formulations of that conjecture.

That can be zetas, for which the Riemann hypothesis is not true, shows
the example of the {\bf Hurwitz zeta function} $\zeta(s,a) , a \in
(0,1]$. It is defined for $Re(s)>1$ as the translated by $a$ Diriclet
series (the generating function of $\LN^{*}\;+\;a$) :
\begin{equation}
\zeta(s,a)\;=\;\sum_{n=0}^{\infty}\frac{1}{(n+a)^{s}}\;,\;0<a \le 1.
\end{equation}
The Hurwitz zeta $\zeta(s,a)$ has the analytic continuation for the
whole complex plane. Let
\begin{displaymath}
\eta\;=\;\frac{\sqrt{10-2\sqrt{5}}\;-\;2}{\sqrt{5}\;-\;1}
\end{displaymath}
and let
\begin{equation}
\zeta_{H}(s)\;:=\;\frac{1}{5^{s}}(\zeta(s,1/5)+\eta \zeta(s,2/5)-\eta
\zeta(s,3/5)-\zeta(s,4/5).
\end{equation}
Then $\zeta_{H}$ is an entine function of $s$ and satisfies the
following functional equation 
\begin{equation}
(\frac{5}{\pi})^{s/2}\Gamma(\frac{1+s}{2})\zeta_{H}(s)\;=\;(\frac{5}{\pi})
^{(1-s)/2}\Gamma(1-\frac{s}{2})\zeta_{H}(1-s).
\end{equation}
Applying some methods from the Riemann zeta function, one can show that
$\zeta_{H}$ has {\bf infinitely many zeros on the critical line} (the
Hardy-Littlewood type theorem for $\zeta_{H}$). However $\zeta_{H}$ has
also {\bf infinitely many zeros in the half-plane $Re(s)>1$}. The
essential difference between $\zeta_{H}$ and $\zeta$ consists on the
fact that $\zeta_{H}$ {\bf does not have the Euler product}! Hence the
corollary, that to have the full Riemann hypothesis, there existence
only a functional equation is not suffices (at least in the area
$Re(s)>1$) and the fact that $\zeta$ has the Euler product is
important. More exactly, the (RH) for $\zeta_{\LQ}$ splits into two different
kind of conjectures: the {\bf Right RH}
\begin{displaymath}
(RH_{1^{+}})\;\zeta_{\LQ}(s) \ne 0\;for\; Re(s)>1
\end{displaymath}
and the {\bf Left RH} 
\begin{displaymath}
(RH_{1^{-}})\;\zeta_{\LQ}(s) \ne 0\;for\;Re(s)\in [0,1].
\end{displaymath}

(6) {\bf The p-adic Riemann hypothesis} (see [WH]).

Researchers working in the theory of zetas are intrigued by the fact
that the Riemann hypothesis arises as a conjecture for the quite
another type of zetas : the zeta functions of divisors, introduced by
{\bf Wan} in [W1] (see also [W2]).

Let $\LF_{q}$ be a finite field with $q$ elements, $q$ a power of a
prime $p$. Let $X$ be a projective n-dimensional integral scheme
defined over $\LF_{q}$. Let $0 \le r \le n$ be an integer. A prime
r-cycle of $X$ is an r-dimensional closed integral subscheme of $X$
defined over $\LF_{q}$. An r-cycle on $X$ is a formal finite linear
combination of prime r-cycles. An r-cycle is called effective, denoted
$\sum n_{i}P_{i}\le 0$, if each $n_{i} \le 0$.

Each prime r-cycle $P$ has an associated graded coordinate ring
$\oplus_{k=0}^{\infty}S_{k}(P)$ since $X$ is projective. By a theorem
of Hilbert-Serre, for all $k$ sufficiently large, we have
$dim_{\LF_{q}}S_{k}(P)$ equal to a polynomial $a_{r}k^{r}+$(lower
terms). Define the degree of $P$, denoted $deg(P)$, as $r!$ times the
leading coefficient $a_{r}$. We extend the definition of degree to
arbitrary r-cycles by $deg(\sum n_{i}P_{i}) := \sum n_{i}deg(P_{i})$.

Defining the degree allows us to measure and compare the prime
r-cycles. Define the zeta function of algebraic r-cycles on $X$ as
\begin{equation}
Z_{r}(X,T)\;:=\;\prod_{P}(1\;-\;T^{deg(P)})^{-1},
\end{equation}
where the product is taken over all prime r-cycles $P$ in $X$.

Denote the set of all effective r-cycles of degree $d$ on $X$ by
$E_{r,d}(X)$. A theorem of Chow and van der Waerder states that this
set has the structure of a projective variety. Since we are over a
finite field, $E_{r,d}(X)$ is finite. This means that $Z_{r}(X,T)$ is a
well-defined element of $1\;+\;T \LZ[[T]]$, and so, converges {\bf
p-adically} in the open unit disk $\mid T \mid_{p}<1$.

Equivalent forms of this zeta function are
\begin{displaymath}
Z_{r}(X,T)=\sum_{d=0}^{\infty}\#E_{r,d}(X)T^{d}=\prod_{d=1}^{\infty}(1-T
^{d})^{-N_{d}}=exp(\sum_{k=1}^{\infty}\frac{T^{k}}{k}W_{k}),
\end{displaymath}
where $N_{d}$ is the number of prime r-cycles of degree $d$ and
$W_{k}:=\sum_{d \mid k}d N_{d}$ is the weighted number of prime r-cycles
of degree dividing $k$, each prime r-cycle of degree $d$ is counted $d$
times.

This zeta function will not be complex analytic in general, unlike the
above classical zeta functions. However, {\bf the conjectural p-adic
meromorphic continuation} immediately yields a p-adic formula in terms
of the p-adic zeros and poles of this zeta function as we shall see
below. 

The p-adic meromorphic continuation of $Z_{r}(X,T)$ would imply the
complete p-adic factorization
\begin{equation}
Z_{r}(X,T)\;=\;\frac{\prod(1\;-\;\alpha_{i}T)}{\prod(1\;-\;\beta_{j}T)},
\end{equation}
where the products are now infinite with $\alpha_{i}\rightarrow 0$ and
$\beta_{j} \rightarrow 0$ in $\LC_{p}$. Taking the logarithmic
derivative, we obtain a formula for $W_{d}$ in terms of infinite
series:
\begin{displaymath}
W_{d}\;=\;\sum \beta_{j}^{d}\;-\;\sum\alpha_{i}^{d}.
\end{displaymath}
Thus, the p-adic meromorphic continuation implies a well-structured
formula for the sequence $W_{k}$. By Mobius inversion, this gives a
well-structured formula for the sequence $N_{k}$ as well.

If $Z_{r}(X,T)$ is p-adic meromorphic, we can adjoint all the
reciprocal zeros $\alpha_{i}'s$ and all the reciprocal poles
$\beta_{j}'s$ to $\LQ_{p}$. The resulting field extension of $\LQ_{p}$
is called the {\bf splitting field} of $Z_{r}(X,T)$ over $\LQ_{p}$.
This splitting field is automatically a Galois extension (possibly of
infinite degree) over $\LQ_{p}$ by the Weierstrass factorization of
$Z_{r}(X,T)$ over $\LQ_{p}$ and the fact that we are in characteristic
zero.

Let us observe that the congruence Weil zeta Z(X,s) is exactly the case
$Z_{0}(X,T)$. Thus $Z_{0}(X,T) \in \LQ(T)$ is rational and it satisfies
the Riemann hypothesis by the Weil Conjectures. If $r = n = dim(X)$,
then $Z_{n}(X,T)$ is trivially rational, its zeros and poles are roots
of unity, and thus satisfies the Riemann hypothesis as well. In
particular, $Z_{r}(X,T)$ is well understood if $n \le 1$. So, we will
assume $n = dim(X) \ge 2$ and $1 \le r \le n-1$.

Let $CH_{r}(X)$ be the Chow group of r-cycles on $X$; that is, the free
abelian group generated by the prime r-cycles on $X$ modulo the
rational equivalence. Let $EffCone_{r}(X)$ be the set of effective
r-cycle classes in $CH_{r}(X)$. Assume that $EffCone_{r}(X)$ is a
finitely generated monoid. Then is postulated the following {\bf p-adic
Riemann hypothesis}: the splitting field of $Z_{r}(X,T)$ over $\LQ_{p}$
is a finite extension of $\LQ_{p}$. Moreover, all zeros and poles,
except for finitely many, are simple.

(8) {\bf Dynamical zeta functions}.

Let $f : X \longrightarrow X$ be such a {\bf homeomorphism} of a {\bf
topological space} X that the number $N_{n}(f)$ of the points of the
period $n$, i.e. the number of solutions of the equation $f^{n}(x) = x$
(obviously $f^{n}$ denotes the n-th iteration of $f$) is finite for
each $n \ge 1$. The function (or rather a formal series)
\begin{equation}
\zeta_{f}(t)\;=\;exp(\sum_{n=1}^{\infty}\frac{N_{n}(f)t^{n}}{n})
\end{equation}
is called the {\bf zeta function of a homeomorphism} $f$ (see [AI]). It
is an extremaly surprising fact, that similarly as in the algebraic
geometry case of congruence Weil zeta , such dynamical zetas are often rational. 
More exactly, let $(\Sigma_{A}, f)$ be an irreducible {\bf Topological
Markov Chain} (TMC in short, see [AI]) with the periodicity index $h$.
Then (see [AI, Th.4.2])
\begin{displaymath}
(1)\;\;N_{n}(f)\;=\;h \lambda^{n}(A)\;+\;\sum_{\mid \lambda_{i} \mid <
\lambda(A)}\lambda_{i}^{n}\;\;for\;\;n=ph,
\end{displaymath}
\begin{displaymath}
(2)\;\;N_{n}(f)\;=\;0\;\;if\;\;n \simeq 0(mod h),
\end{displaymath}
The zeta function $\zeta_{f}(t)$ of an arbitrary TMC $(\Sigma_{A}, f)$ is 
rational and the Riemann hypothesis is trivial in his "zero part" :
\begin{displaymath}
\zeta_{f}(t)\;=\;\frac{1}{det(E\;-\;tA)},
\end{displaymath}
and the defining series is convergent for $\mid t \mid <1/\lambda(A)$
(here $\lambda(A)$ denotes the largest positive eigenvalue of the
matrix $A$).

{\bf P. Bowen} has also showed that zeta functions of $A^{\#}$-homeomorphism
are rational (see Th.12.1 and Th.13.1 of [AI]). Extremaly exciting is also
an application of dynamical zetas for {\bf Lorentz attractors} (see [W]).
They appear in the proof of the fundamental {\bf Guckenheimer theorem}, which
explains the structure of the Lorentz attractors (see[W, Theorem]) : there
exists only two different topological types of the Lorentz attractors.

More about the dynamical-topological zeta functions a reader can find in a 
beautiful reviewed article [Ru] by D. Ruelle.

Let us observe that dynamical zetas are strictly connected with congruence
Weil zetas : let $Fr_{q} :V \longrightarrow V$ be the {\bf Frobenius map}
(acting by $z \longrightarrow z^{q}$ on coordinates), $\mid Fix f^{n} \mid$
is the number of fixed points of the n-th iterate of $f$. Then
\begin{displaymath}
Z(V/\LF_{q})(s)
\;=\;\zeta_{Fr_{q}}(s)\;=\;exp(\sum_{n=1}^{\infty}\frac{s^{n}}{n}\mid
Fix Fr_{q}^{n} \mid).
\end{displaymath}

Finally, the last interesting future of the dynamical zetas $\zeta_{f}(s)$
is the {\bf triviality in zeros} (in general) of the Riemann hypothesis in 
this case (although they have "non-trivial" functional equations and Euler 
products in this case). However in the much general zetas - like zeta functions
associated with the weighted dynamical systems, for example the {\bf 
Lefschets zeta function} $\zeta_{L}(s)$, the {\bf Selberg zeta function}
$\zeta_{S}(s):=\prod_{k=0}^{\infty}\zeta(s+k)^{-1}$ and the currently
popular {\bf Ihara-Selberg zeta function} $\zeta_{I}(s)$ associated
with a {\bf finite unoriented graph} $G$ - the Riemann hypothesis problem
is actual (well-posed) (see [Ru]). It is also well-known that $1/\zeta_{I}$
is a {\bf polynomial} and that $\zeta_{I}$ satisfies the {\bf Riemann
hypothesis precisely when $G$ is Ramanujan} ( Ramanujan graphs were
named by Lubotzky, Phillips, and Sarnak; examples are not easy to
construct).

In [MA], [MH], [AM] and [ML] we initiated and developed the method (technique)
of the proving of the {\bf Riemann Hypothesis}( RH in short) by using some 
methods from the {\bf measure and integration theory} both on LCA groups (see 
[MA] and [MH]) and infinite-dimensional linear topological spaces (see [AM] 
and [MR]).

The mentioned above meromorphic extension $\zeta(s)$ of the Dirichlet
series (1.1) is given explicite by the following classical {\bf Riemann
functional analytic continuation equation}(Rface in short)
\begin{equation}
\zeta(s)\;=\;\frac{\pi^{s/2}}{\int_{0}^{\infty}s^{s/2-1}e^{-x}dx}
[\frac{1}{s(s-1)}\;+\;\int_{1}^{\infty}(x^{s/2-1}+x^{-(s+1)/2})
(\sum_{n=1}^{\infty}e^{-\pi n^{2}x})dx].
\end{equation}

Based on (1.10) and the existence of the {\bf Hodge measure} $H_{2}^{*}$
in [ML] we give a {\bf Hodgian proof} of (RH).

The whole our approach to (RH) (also in this paper) is based on a
subsequence extensions of that principial (Rface) to the more general
class of functions than the {\bf gaussian canonical function} $G(x) :=
e^{-\pi x^{2}}$.

In [MA] we proposed to extend (1.10) and consider the class of {\bf fixed
points} $\omega^{+}$ of the {\bf canonical Fourier transform} ${\cal
F}(f)(x) := \int_{\LR}e^{2 \pi ixy}f(y)dy$ (let us observe that $G$ is
only one special example 0f $\omega^{+}$) of the form (ffp-equation):
\begin{equation}
(M(\omega^{+})\zeta)(s)\;=\;\frac{\omega^{+}(0)}{s(s-1)}\;+\;\int_{1}
^{\infty}\theta(\omega^{+})(x)(x^{s-1}+x^{-s})dx ,
\end{equation}
where $M(f)(s):=\int_{0}^{\infty}x^{s-1}f(x)dx$ stands for the {\bf
Mellin transform} and $\theta(f)(x):=\sum_{n=1}^{\infty}f(nx)$ is the
{\bf theta Jacobi transform}.
Beside (1.11), the second important technical tool is the notion of the
{\bf RH-fixed point of} {\cal F} : $\omega^{+}_{A} - G = A$ (associated
with an amplitude $A$) (see [MA] for details).

To prove the existence of $\omega^{+}$ we used an expanded aparatus of
the measure theory - in particular {\bf Haar measures}, {\bf Riesz
measures}, {\bf Bogoluboff-Kriloff measures} and {\bf Herbrandt
distributions} - to give an {\bf algebraic proof} of (RH).

That direction of investigations was next directly continued in [MH].
The ffp-equation (1.11) is subsequently extended for classes of {\bf
eigenvectors} of the whole {\bf family} $({\cal F}^{p}:p>0)$ of {\bf
Fourier transforms} $({\cal F}^{p}f)(x) = \int_{\LR}e^{2 p^{2}
ixy}f(y)dy$, corresponding to the {\bf real valued eigenvalues}
$\lambda : {\cal F}^{p}\omega_{\lambda}^{p} = \lambda \omega_{\lambda}^{p}$.
Thus, in [MH] we derived the ${\cal F}$-eigenvalue analytic continuation 
equations for {\bf eigenvectors} $\omega_{\lambda}^{p}$ of ${\cal
F}^{p}$, i.e. ${\cal F}^{p}\omega_{\lambda}^{p} = \lambda
\omega_{\lambda}^{p}$ ( although for the purposes of [MH] we need only
some very special eigenvectors - the {\bf Hermite functions}
$H_{p}^{n}(x)$ defined simply by the formula : let $G_{p}(x) := e^{-p^{2}x^{2}}$ 
be the gaussian function with the parameter $p>0$. Then $H_{p}^{n}(x)
:= G_{p}^{(n)}(x), n=0 , 1, 2, ... ,$ are eigenvectors of ${\cal F}^{p}$
corresponding to the eigenvalues $\lambda = \pm \frac{\sqrt{\pi}}{p}$.
Thus we have : let us denote $m(s)=1-s$. Then
\begin{equation}
\sum_{g \in \{id_{\LC}, m\}}\mid \lambda
\mid^{-2g(s)}(M(\omega^{p}_{\lambda})\zeta)(g(s))\;=\;\frac{(\lambda+1)
\omega_{\lambda}^{p}(0)}{2s(s-1)}\;+\;
\end{equation}
\begin{displaymath}
(\mid \lambda
\mid^{-1}\;+\;sgn(\lambda))\int_{1}^{\infty}\theta(\omega_{\lambda}^{p})(
\lambda^{2}x)(x^{s-1}\;+\;sgn(\lambda)x^{-s})dx.
\end{displaymath}
In some sense, in this paper we consider a "maximal" extension of (1.10)
- the so called {\bf Muntz's relations}- for the functions $p$ from the
cylinder of the {\bf Poisson space} $\LP_{1}$ (see Th.2), of the form
\begin{equation}
(M(p)\zeta)(s)\;=\;\frac{1}{s(s-1)}+\int_{1}^{\infty}[x^{s-1}\theta(p)(x)+
x^{-s}\theta(({\cal F}p))(x)] dx.
\end{equation}

A big Poisson space $\LP_{1}(\LR_{+})$ has such an adventage that the
famous {\bf Poisson Summation Formula} (PSF in short) holds for the functions
from that one. Moreover, on the Borel $\sigma$-field ${\cal P}$ of
$\LP_{1}$, we will be able to define - the principial for this short
proof of (RH) - the family of the {\bf Wiener-Riemann measures} $\{r_{m}\}$ -
indexed by some simple class of moment functions $m$.

The family $\{r_{m}\}$ is in fact induced by the {\bf standard Brownian
motion} $B^{0}$ (see Section 2) - and hence follows the huge role of
the Brownian motion (or rather the {\bf Wiener measure} $w$) for that
proof of (RH).

In the above three types of the functional equations written (mentioned)
above, appears always the simple algebraic meromorphic function
\begin{displaymath}
\zeta_{\LP^{1}(\LC)}(s)\;:=\;\frac{1}{s(s-1)},
\end{displaymath}
(the congruence Weil zeta of the complex projective space $\LP^{1}(\LC)$
- see the end of Section 3).

Subsequently, in $\zeta_{\LP^{1}(\LC)}(s)$ is written the fundamental
polynomial of two variables with integer coefficients and degree 2 :
\begin{equation}
\zeta_{t}(s)\;=\;Im(\frac{1}{s(s-1)})\;=\;Im(s)(2Re(s)\;-\;1)\;=\;y(2x\;-\;1).
\end{equation}

According to the below surprising property of $\zeta_{t}$ and for some
another important algebraic reason, we call the function $\zeta_{t}(s)$
- the {\bf trivial zeta} - since its formally satisfies the following
{\bf Trivial Riemann Hypothesis} (TRH in short) - according to its
formal reasonblance (similarity) with the true (RH).
\begin{displaymath}
(TRH)\;\;(\zeta_{t}(s)=0)\land (Im(s) \ne 0) \Longrightarrow Re(s)=1/2.
\end{displaymath}

As a consequence, in [MA] we formulated the following {\bf Main
Algebraic Conjecture} (MAC in short) that :

{\bf (MAC) $\;\;TRH\;\;implies\;\;RH$}.

More exactly, we can formulate (MAC) in the following algebraic geometry
language : let $\zeta(\LC) := \{s \in \LC : \zeta(s) = 0\}$ be the {\bf
zero-dimensional comlex analytic manifold} and $\zeta_{t}(\LC) :=
\{(x,y) \in \LR^{2} : \zeta_{t}(x,y) = 0\}$ be the {\bf 1-dimensional
algebraic (affine) variete over $\LR$}. Then (MAC) means that
\begin{equation}
(MAC)\;\;\;\zeta(\LC)\;\subset \zeta_{t}(\LC),
\end{equation}
i.e. $\zeta(\LC)$ is a {\bf subvariete} of $\zeta_{t}(\LC)$.

The Muntz's relations considered in this paper, are modeled on the
following extension of the (Rface) considered in [AM] : let $\theta MA$
be the $\LR_{+}$-cone of {\bf theta-Mellin admissible} real valued
functions on $\LR_{+}$, i.e. the smallest space of functions on which
the Mellin transform $M$ {\bf does not vanish} and the (PSF) holds.
Moreover, each $f \in \theta MA$ satisfies the initial condition :
$\hat{f}(0) = f(0) = 0$. As a consequence of the above conditions on
$\theta MA$ we get the following family of (Rface) indexed by $\theta
MA$ :
\begin{equation}
(M(f)\zeta)(s)\;=\;\int_{1}^{\infty}(x^{s-1}\theta(f)(x)\;+\;x^{-s}
\theta(\hat{f}(x))dx ,
\end{equation}
i.e. in (1.16) does not appear the congruence zeta $\frac{1}{s(s-1)}$.

In opposite to (MAC) we thus also have the {\bf Main Oscilatory
Conjecture} (MOC in short) : in the right-hand sides of (1.10), (1.11),
(1.12) and (1.13) appears another very importantant {\bf analytic 
expressions} - the {\bf Fresnel integrals} $F(A)(\nu) := \int_{0}^{\infty}sin \nu x
dA(x)$, associated with a {\bf positive measure} $A$ and a {\bf
periodicity} $\nu >0$. {\bf Arnold} and others devoted the whole book
[AWG-Z] to the asymptotics of the so called {\bf osccilatory integrals}
of the form (see [AWG-Z, Sect.2 and 3])
\begin{displaymath}
\int_{\LR^{n}}e^{i(\phi(x))}A(x)dx,
\end{displaymath}
where the functions $\phi$ and $A$ are called the {\bf phase} and {\bf 
amplitude}, respectively.

The main property of $F(A)(\nu)$ is its {\bf positivity} ( the Fresnel lemmma
- Oscilatory lemma or finally - the Nakayama lemmma) - for the {\bf positive, 
continuous, integrable} and {\bf decreasing} densities $(\frac{dA}{dx}$.

The (MOC) is also a simple logical statement that :

{\bf The Fresnel lemmma} $\;\;$implies$\;\;$RH.

(It is mainly proved in [AM]).

Thus, the whole effort of the papers [MA], [MH] and [ML] was
concentated on the proving of (MAC)+(MOC). Moreover - in fact - we only
proved and used (MOC) in [AM] and only (MAC) - in this paper.

In the all mentioned papers devoted to the proofs of (RH), the proof
was always attained by the writting of a suitable {\bf Riemann hypothesis
functional equations}. Exactly the same is in this paper, the Riemann 
hypothesis is attained by writing the following {\bf Brownian motion
Riemann hypothesis equation}:
\begin{equation}
\lim_{k}Im(m_{s}(\LP_{1},{\cal P},
r_{m_{k}})\zeta(s))\;=\;\frac{\zeta_{t}(s)}{\mid \zeta_{\LP^{1}(\LC)}(s) \mid
^{2}},
\end{equation}
with $Re(s) \in (0,1/2)$.

In such a way, the {\bf measure characteristic numbers} $m_{s}$ (see Sect.4)
gives the direct relation beetwen the {\bf analytic properties} of $\zeta$ and 
the {\bf algebraic nature} of the pair $(\zeta_{t}, \zeta_{\LP^{1}(\LC)})$. One
can say that in the family $(m_{s}, b_{s})$ (see Sect.2) is written the 
{\bf whole information} concerning the {\bf Riemann hypothesis} - like in the {\bf
Betti numbers} $B^{i}(X) := dim_{\LC} H^{i}(X,\LC)$ ( where
$H^{i}(X,\LC)$ is the i-th cohomology space of a manifold $X$ with a
coefficients in $\LC, i = 0, 1, ... ,2n ; n = dim(X)$)  and in the
{\bf Euler-Poincare characteristic} $\chi(X) :=
\sum_{i=0}^{2n}(-1)^{i}B^{i}(X)$ (i.e. in the pair $(B(X),\chi(X))$) is 
written the {\bf whole information} concerning the {\bf congruence Riemann 
hypothesis} (see ($W_{R}$) and (cRH)).

\section{The Poisson space $\LP_{1}(\LR_{+})$ and the Wiener-Riemann
measures $r_{m}$.}
Let $C = C(\LR_{+})$ be the real vector space of all real valued {\bf
continuous} functions on $\LR_{+}=[0,+\infty)$. Let $L^{1} =
L^{1}(\LR_{+}, dx)$ be the real Banach space of all absolutely integrable
real valued functions. If $f$ belongs to $C \cup L^{1}$ then by $f^{+}$
we denote its {\bf symmetrisation}, i.e. 
\begin{equation}
f^{+}(x)\;:=\;\{\frac{f(x)\;if\;x\ge 0}{f(-x)\;if\;x<0}.
\end{equation}
By ${\cal F}$ we denote the {\bf canonical Fourier (cosine) transform}, i.e. 
\begin{equation}
({\cal F}f^{+})(x)\;:=\;\int_{-\infty}^{+\infty}e^{2\pi
ixy}f^{+}(y)dy\;=:\;\hat{f^{+}}(x)=2\int_{0}^{\infty}cos(2\pi xy)f(y)dy;
\end{equation}
$x\in \LR, f \in L^{1}$.

Moreover, we denote
\begin{equation}
\hat{C}\;:=\;\{f : \LR_{+}\longrightarrow \LR : \hat{f^{+}}\mid \LR_{+}
\in C \},
\end{equation}
and
\begin{equation}
\hat{L^{1}}\;:=\;\{f :\LR_{+}\longrightarrow \LR : \hat{f^{+}}\mid
\LR_{+} \in \LR_{+}\}.
\end{equation}
Finally, we denote the real vector space being the {\bf intersection}
of the all above Banach and Frechet spaces by $Inv(\LR_{+})$, i.e.
\begin{equation}
Inv(\LR_{+})\;:=\;C \cap \hat{C} \cap L^{1} \cap \hat{L}^{1}.
\end{equation}

As in the all previous our papers devoted to (gRH) (see [MA,MH,AM,ML]
and [MR]) , it will be very convenient to consider the {\bf canonical
Jacobi theta transform} $\theta$.

It is very convenient to consider $\theta$ on the {\bf Schwartz space}
${\cal S}(\LR_{+})$ of all rapidly decreasing and smooth functions on
$\LR_{+}$. Thus
\begin{equation}
\theta(f)(x)\;:=\;\int_{\LN^{*}}f(nx)dc_{\LZ}(n)\;=\;\sum_{n=1}^{\infty}f
(nx),
\end{equation}
where $f \in {\cal S}(\LR_{+}), x>0$  and $c_{\LZ}$ denotes the {\bf
calculating measure} on $\LN^{*}$, i.e. the unique Haar measure on the
LCA group $(\LZ, +)$ normalized by $c_{\LZ}(\{0\})=1$.

As it was communicated to us by {\bf S. Albeverio} (by a private
communication), the Schwartz space is to big for the below {\bf Mellin
transform} $M$. It follows from the fact that $M$ has also a
"singularity" at zero (beside the singularity at the plus infinity).

Thus by ${\cal M}(\LR_{+})$ we denote the sub-vector space of ${\cal
S}(\LR_{+})$ consisting with such $f \in {\cal S}(\LR_{+})$ that
\begin{displaymath}
(MT)\;\;\int_{0}^{1}x^{Re(s)-1}\mid f(x)\mid dx\;<\;+\infty,
\end{displaymath}
for all $Re(s) \in I:=(0,1)$.
\begin{re}
Let us observe that ${\cal M}(\LR_{+})$ is essentially smaller than
${\cal S}(\LR_{+})$. Reely, for example the n-th {\bf Hermite functions} 
functions $H_{n}(x)$ being the n-th derivatives of the canonical
gaussian functions $G(x)$, i.e. $H_{n}(x):= G^{(n)}(x)$ obviously belong to 
the Schwartz space. But their Mellin transform
\begin{displaymath}
M(H_{n})(s):=\int_{0}^{\infty}x^{s-1}H_{n}(x)dx=\int_{0}^{\infty}x^{s-1}
G^{(n)}(x)dx=s(s-1) ... (s-n)\Gamma(s-n)
\end{displaymath}
is {\bf convergent} for $s$ with $Re(s)  > n$!(see also [MH, Prop.7])

Thus we have deal with rather subtle problem of the convergence of
non-proper integrals. They also have the "singularity" at zero (besides
at infinity) and obviously the Schwartz space ${\cal S}(\LR)$ is a
"good" condition on the behaviour of functions in infinity.
\end{re}

The space ${\cal M}(\LR_{+})$ we call the space of {\bf Mellin
admissible functions}. The {\bf Mellin transform} $M : {\cal
M}(\LR_{+})\longrightarrow \LC$  is defined as
\begin{equation}
M(f)(s)\;:=\;\int_{0}^{+\infty}x^{s}f(x)\frac{dx}{x}\;;\;Re(s)>0, f \in
{\cal M}(\LR_{+}).
\end{equation}

We propose the following fundamental for the purposes of this paper 
\begin{de}
By the {\bf Poisson space} $\LP_{1}(\LR_{+})$ ( $\LP_{1}$ for short) we
mean (understand) the set of all functions $p: \LR_{+}\longrightarrow
\LR$ which satisfies the following five conditions :

($P_{0}$) $p(0) = 1$,

($P_{1}$) $p^{+} \in Inv(\LR)$,

($P_{2}$) for all $x, y \in \LR_{+}$
\begin{displaymath}
\sum_{n=1}^{\infty}x^{2}\mid p(x+ny)
\mid\;=:\;\theta(x^{2}p^{+}_{+x})(y)<\infty,
\end{displaymath}
($P_{3}$) for all $x, y \in \LR_{+}$
\begin{displaymath}
\sum_{n=1}^{\infty}x^{2}\mid \hat{p}^{+}(x+ny)
\mid\;=:\;\theta(x^{2}\hat{p}^{+}_{+x})(y)<\infty,
\end{displaymath}
($P_{4}$) $p \in {\cal M}(\LR_{+})$.
\end{de}

The importance of the Poisson space $\LP_{1}$ follows from the two
sources : 

(1) for the functions from $\LP_{1}$ hold the below {\bf Poisson Summation 
Formula}( fundamental when we work with zetas) (PSF in short) and the Mellin 
transform is well-defined on its.

(2) On $\LP_{1}$ there exist probabilities $r_{m}$ - with properties,
which are sufficient to prove the Riemann Hypothesis. More exactly, we
have the following technical proposition and theorem :
\begin{pr}
For each $p \in \LP_{1}$ holds

(i) the {\bf Poisson Summation Formula} (PSF in short) , i.e.
\begin{displaymath}
(PSF)\;\;\theta(p)(1)+p(0)\;=\;\theta(\hat{p^{+}})(1)+\hat{p^{+}}(0),
\end{displaymath}
and

(ii) $M(p)(s)$ is well-defined (exists) for all $s$ with $Re(s)>0$ and
moreover if $p = exp^{-1} \in \LP_{1}$ then
\begin{displaymath}
M(exp^{-1})(s)\;=\;\int_{0}^{\infty}x^{s-1}e^{-x}dx\;=\;\Gamma(s),
\end{displaymath}
is the classical gamma function which {\bf does not vanish everywhere}.
\end{pr}
All the above notions and facts, as well as their proofs, in a large context 
and deepnest can be find in the Weil's book [We].

For the purposes of this paper we fix an arbitrary {\bf continuous
moment function} $m : \LR_{+} \longrightarrow \LR_{+}$, which satisfies
the following three conditions 

($m_{0}$) $m(0) = 1$,

($m_{1}$) $0 \le m(x) \le 1$ for all $x \in \LR_{+}$

($m_{2}$) $m(x)=0$  for $x \ge 1$.

A {\bf stochastic process} $B^{m} = (B_{t}^{m}: t \ge 0)$ defined on a 
probability space $(\Omega, {\cal A}, P)$ is said to be the {\bf Brownian 
motion with a moment function $m$} if it is {\bf gaussian} (i.e. its
all 1-dimensional distributions are gaussian) , its medium function is
$m$ :
\begin{equation}
EB_{t}^{m}\;=\;m(t),
\end{equation}
and its {\bf covariance function} (or correlation function)
\begin{equation}
E(B_{t}^{m}-m(t))(B_{s}^{m}-m(s))\;=\;min(t,s).
\end{equation}
Here and all in the sequel $EX$ denotes the {\bf expected value} of a
random variable $X$ (rv in short).

For there existence of such $B^{m}$ see e.g. [Wong, II.3].
In particular, $B^{m}$ has got the {\bf continuous paths}
(trajectories), since it satisfies the well-known {\bf Kolmogorov
condition} (see e.g. [Wong, II.4, Th.4.2]).

Let us observe at once that $B^{m}$ has the direct representation : let
$B^{0}=(B^{0}_{t}:t\ge 0)$ be the {\bf standard Brownian motion}. Then
\begin{equation}
B^{m}_{t}\;=\;B_{t}^{0}\;+\;m(t).
\end{equation}

The Brownian motion $B^{m}$ can be considered as a {\bf random element}
(re in short) $B^{m} : (\Omega, {\cal A}, P)\longrightarrow (C, {\cal C})$
, where ${\cal C}$ is the {\bf $\sigma$-field generated by cylinders in
$C$}, .i.e. the sets of the form
\begin{displaymath}
C\;=\;C(t_{1}, ... , t_{n};A)\;:=\;\{c \in C : (c(t_{1}), ...,c(t_{n})) \in
A\},
\end{displaymath}
where A is a Borel set from $\LR^{n}$ and $0<t_{1}< ... <t_{n}$.

It is well-known that ${\cal :C}$  coincides with the $\sigma$-field of
all Borel sets ${\cal B}$ of $C$ endowed with the Frechet topology of
the almost uniform convergence on compacta. Normally, $(C, {\cal C})$
is called the {\bf phase space} of $B^{m}$.

Finally, $B^{m}$ (as a re) determines the {\bf probability}  $w^{m}$ on
the phase space $(C, {\cal C})$ by the formula (the distribution of
$B^{m}$) :
\begin{equation}
w^{m}(C)\;:=\;P((B^{m})^{-1}(C)) \;\;,\;C\in {\cal C}.
\end{equation}
In the sequel we call $w^{m}$ - the {\bf m-Wiener measure}. The
m-Wiener measure $w^{m}$ (similarly like p-stable measures $s^{p}$ with
$p \in (0,1)$ considered in [AM]) permits us to define principial for
this short proof of (RH) - the probability $r_{m}$.

To do this, we first of all define the sequence of the {\bf Wiener-Riemann 
processes} $\{R^{n} : n=0,1, ...\}$ by the simple formula
\begin{equation}
R^{n}_{t}\;:=\;G(t)B^{m}_{\sqrt{t}},\;\;t \in [n,n+1].
\end{equation}

Let ${\cal P}$ be the $\sigma$-field of $\LP_{1}$ generated by cylinders of
$\LP_{1}$.
\begin{th}
For each moment function $m$ which satisfies the conditions $(m_{0})-(m_{2})$, 
there exists a $\sigma$-additive measure $r_{m}$ on the phase space
$(\LP_{1}(\LR_{+}), {\cal P})$ with the following four - let us say -
RH properties):

($RH_{0}$)(Non-triviality, probabilisticity):  $r_{m}(\LP_{1}) = 1$.

($RH_{1}$)(Starting
point):$\int_{\LP_{1}}p(0)dr_{m}(p)\;=\;EB_{0}^{m}\;=\;m(0)\;=1$.

($RH_{2}$)(Vanishing of moments): for all $t \ge 1$ holds
\begin{displaymath}
\int_{\LP_{1}}p(t)dr_{m}(p)\;=\;E(G(t)B_{\sqrt{t}}^{m})\;=0.
\end{displaymath}
($RH_{3}$)({\bf The Fubini obstacle - Hardy-Littlewood theorem obstacle -
Existence of moments of 1/2-stable Levy distributions}):

the {\bf Brownian characteristic numbers $b_{s}$ of $\LP_{1}$ }or {\bf double
Wiener integrals of} $\LP_{1}$)
\begin{equation}
b_{s}=b_{s}(\LP_{1},{\cal P},r_{m};\LR_{+})\;:=\;\int \int_{\LP_{1}\times 
\LR_{+}}\mid x^{-s}p(x)\mid dr_{m}(p)dx
\end{equation}
are

(i) {\bf finite (convergent)}  if $Re(s) \in (0,1/2)$, and

(ii){\bf infinite (divergent)} if $Re(s) \ge 1/2$.

In the sequel the measure $r_{m}$ we call the {\bf Wiener-Riemann
measure} (by analogy with the Levy-Riemann measures $\{r_{p} : p \in
I:=(0,1)\}$ considered in [AM]).
\end{th}
{\bf Proof}. We define the measure $r_{m}$ by the following simple
formula :
\begin{displaymath}
(WRm)\;\;r_{m}(A)\;:=\;\sum_{n=1}^{\infty}\frac{1}{2^{n}}P(R_{(\cdot)}^{n}
\in A \cap C[n-1,n])\;\;,\;\;A \in {\cal P},
\end{displaymath}
where by $C[n-1,n]$ we denoted the space of all real valued continuous
functions defined on the segment $[n-1,n]$, considered as a subspace of
$C$ through the embeding $e_{n} : C[n-1,n] \longrightarrow C(\LR_{+})$
by the formula : $e_{n}(f)(t) = f(t)$ if $n-1\le t \le n;
e_{n}(f)(t)=f(n-1)$ if $t \le n-1$ and $e_{n}(f)(t) =f(n)$ if $t \ge
n$ and moreover we have
\begin{displaymath}
P(R_{(.)}^{n}\in A \cap C[n-1,n])=P(f \in C(\LR_{+}):
G(t)(f(\sqrt{t})=m(t) \in A \cap C[n-1,n])\;=
\end{displaymath}
\begin{displaymath}
\;=\;P(f \in C(\LR_{+}): f(t) \in (A \cap C[\sqrt{(n-1)},\sqrt{n}]-m)G^{-1})
\;=\;w^{0}((A \cap C[\sqrt{n-1}, \sqrt{n}]-m)G^{-1}).
\end{displaymath}

Now the conditions $(RH_{0}-RH_{3})$ are trivially satisfy by $r_{m}$.
Reely, the conditions $(P_{0}-P_{4})$ on the Poisson space ( as the
conditions on $p$ at zero and infinity) are quite "independent" from
the condition $p \in C[n-1,n]$ for $n \ge 2$, i.e. $\LP_{1} \cap
C[n-1,n] = C[n-1,n]$ and $w_{0}(C[n-1,n])=1$.

Moreover
\begin{displaymath}
\int_{\LP_{1}}p(t)dr_{m}(p)\;=\;\frac{1}{2^{n}}\int_{C[n-1,n]}p(t)dw_{m}(p)
=G(t)m(\sqrt{t})/2^{n}=0,
\end{displaymath}
if $t \in [n-1,n], n \ge 1$.

Thus, we concentrate only on the unique non-trivial property $(RH_{4})$
of $r_{m}$. The idea of the proof is taken from (based on) the Proposition 3
of [AM]. In fact, in the Brownian characteristic numbers $b_{s}$ are
written all what is important for our Brownian proof of (RH).

We first give the upper approximation of the iterated integral
$I_{dxdr}(s):= \int_{0}^{\infty}dx(\int_{\LP_{1}})\mid x^{s-1}p(x) \mid
dr_{m}(p)$ for $s$ with $u=Re(s) \in (0,1/2)$. 

Let us observe that
\begin{displaymath}
\int_{\LP_{1}}\mid p(x)\mid dr_{m}(p)\;=\;E(G(x) \mid
B^{m}_{\sqrt{x}}\mid)).
\end{displaymath}
Let $u = Re(s)$. Hence
\begin{equation}
\int_{0}^{\infty}x^{u-1}(\int_{\LP_{1}}\mid p(x) \mid dr_{m}(p))dx\;\le\;
\int_{0}^{\infty}x^{u-1}G(x)m(x)dx
\;+\;\int_{0}^{\infty}x^{u-1}G(x)E\mid B^{0}_{\sqrt{x}} \mid dx .
\end{equation}
Since in the above inequality the first integral obviously exists, thus
the problem of the {\bf convergence} of $I_{dxdr}(s)$ is reduced to the
convergence of the integral 
\begin{equation}
\int_{0}^{\infty}x^{u-1}G(x)E\mid B^{0}_{\sqrt{x}} \mid
dx\;=\;\int_{0}^{\infty}x^{u-1}G(x)(\int_{\LR} \mid y+y_{0} \mid
\frac{e^{-\frac{(y+y_{0})^{2}}{2x}}}{\sqrt{2\pi x}}dy)dx,
\end{equation}
for any $y_{0}>0$, according to the facts that the distribution of
$B^{0}_{\sqrt{x}}$  is gaussian with mean zero and variance $\sqrt{x}$
and the translation invariance of the Lebesgue measure.

Now, let us observe ( what was firstly observed in [AM]), that the
iterated integral in the right-hand side of (2.32) - according to the
approximation : $e^{-y^{2}/2x}e^{-2\mid y \mid y_{0}/2x} \le e^{-\mid
 y \mid y_{0}/x}$  and a suitable substitution can be approximated as follows:
\begin{equation}
\le (\int_{0}^{\infty}x^{u}\frac{e^{-y_{0}^{2}/2x}}{\sqrt{2
\pi}x^{3/2}}dx)(y_{0}^{-2}max_{x\ge0}(x^{2}G(x))\int_{\LR}(\mid t \mid
 e^{-\mid t \mid}dt+max_{x \ge 0}(xG(x))\int_{\LR}e^{-\mid t \mid})dt).
\end{equation}
But now, the function $ d_{y_{0}}(x):=
\frac{y_{0}e^{-y_{0}^{2}/2x}}{\sqrt{2 \pi}x^{3/2}}$  is exactly the
{\bf density} of a random variable $L_{y_{0}}$  with the {\bf
1/2-stable Levy distibution} (with a parameter $y_{0}$). It is
well-known that the distribution of $L_{y_{0}}$ is concentrated on
$\LR_{+}$ and that it is the unique p-stable distribution with $p \in
(0,1)$, which has an elementary analytic and simple formula for density
of power-exponential form (see e.g. also [AM]).

The most important fact concerning 1/2-stable Levy distributions (i.e.
some very specific probability measures on $\LR_{+}$) what we use for
this proof, is the problem of the {\bf existence of moments} of $L_{y_{0}}$.
It is well-known that (see [AM]):

\begin{equation}
(i)\;\;E(L_{y_{0}}^{u})\;=\;\int_{0}^{\infty}x^{u}\frac{y_{0}e^{-y_{0}^{2}
/2x}dx}{\sqrt{2\pi}x^{3/2}}\;<\;+\infty \;if\;0<u<1/2,
\end{equation}
and

(ii) $E(L_{y_{0}})^{u}\;=\;+\infty$ if $u\ge 1/2$.

Combining (3.33) with (3.34) we finally obtain that the iterated integral
$I_{dxdr}(s)$ is {\bf finite} if $Re(s) \in (0,1/2)$. Since obviously
the measures $r_{0}$ and $dx$ are {\bf $\sigma$-finite}, then
according to the {\bf Tonelli theorem}  (see e.g. [M, Th.XIII.18.4]),
the Brownian characteristic numbers $b_{s}$  are {\bf finite} in this
case.

We refer a reader to [M, Sect. XIII] for a very detailed and deep
discussion of the Fubini theorem theory.

It is a surprising fact as a huge role plays Fubini theorem (FT in
short) in this Brownian motion approach to (RH) as well as in the
p-stable strategy of the proof of (RH) in [AM]. With such a huge role
of (FT) we first have met when we had worked over the problem of there
existence of {\bf Sazonov topologies} in [Sa].

Looking only at the above upper approximations (evaluations) of
$I_{dxdr}(s)$ for $Re(s) \in (0,1)$, it is not clear why
$I_{dxdr}(s)=+\infty$ for $Re(s)\ge 1/2$. Even worse, at first glance it
seems that also by that method, we can obtain the finitness of
$I_{dxdr}(s)$ in the mentioned above halfplane.

The below lower approximation - also modeled on the previous one -
however shows that (FT) is {\bf violated} in the case of the triplet :

$(\int \int _{\LP_{1}\times \LR_{+}}\mid x^{s-1}p(x) \mid dr_{0}(p)dx,
I_{dxdr}(s), I_{drdx}(s))$ for $Re(s)\ge 1/2$.

Thus, the {\bf Fubini-Tonelli theorem} stands for (is) a true barrier
(obstacle) for the extension of the below $(Rhfe_{B})$ for $Re(s)\ge 1/2$
and is mainly responsible for the {\bf non-triviality} of the {\bf
Riemann hypothesis}.

It is easy to see that (let assume that $y_{0}>1$)
\begin{equation}
b_{s}\;=\; \int \int_{\LP_{1}\times \LR_{+}}\mid x^{s-1}p(x) \mid
dr_{m}(p)dx\; \ge\; \int_{1}^{\infty}x^{u-1}G(x)E\mid B^{0}_{\sqrt{x}} \mid
dx \; \ge\;
\end{equation}
\begin{displaymath}
\;\ge\;(min_{x\ge 1}(x^{2}G(x))\int_{\LR}\mid t \mid e^{-\mid t
\mid}dt\;+\;min_{x \ge 1}(xG(x))\int_{\LR}e^{-\mid t \mid}dt)E(L_{y_{0}}^{u}). 
\end{displaymath}
Thus really $I_{dxrd}(s) = +\infty$ for $u=Re(s) \ge 1/2$, since for
such $s$ there is : $E(L_{y_{0}}^{u}) = \infty$ (see above (ii)).

\section{The Muntz relations for the quintet $(\zeta, \zeta_{\LP^{1}(\LC)}, M, {\cal F}, \theta)$.}

{\bf Muntz} (see [Ti, p.]) was probably the first, who considered the
below extension of the classical {\bf Riemann functional continuation
equation}, for a large class of functions than $exp^{-1}$.

The importance and a huge role of this extension was observed
subsequently in [MA], [MH], [MD] and [AM] - in the different contexts,
for different class of functions. It also seems that it can be compared
only with the ingenuous idea of {\bf Grothendieck} - during his work on
the {\bf congruence Riemann hypothesis} - of the extension of the notion
of the set topology to the category topologies and its main consequence
- there existence of good Weil cohomologies - e.g. l-adic etale
Grothendieck cohomologies.
\begin{th}({\bf Family of Riemann functional analytic continuation
equations for $\zeta$ and the class $\LP_{1}$}).

For each $p \in \LP_{1}$ and $s$ with $Re(s)>0$ the functional equation holds
\begin{equation}
(M(p)\zeta)(s)\;=\;\frac{1}{s(s-1)}\;+\;\int_{1}^{\infty}[x^{s-1}\theta(p)
(x)\;+\;x^{-s}({\cal F}p)(x)]dx.
\end{equation}
\end{th}
{\bf Proof}. Let $p \in \LP_{1}(\LR_{+})$ be arbitrary. By ($P_{5}$), the
Mellin transform $M(p)(s)$ is well-defined for $Re(s)>0$. From the definition
of the Mellin transform $M$ as the integral, on substituting $nx$ for $x$
under the integral, we have
\begin{equation}
\frac{M(p)(s)}{n^{s}}\;=\;\int_{0}^{\infty}p(nx)x^{s-1}dx\;,\;Re(s)>0.
\end{equation}
Hence, for $Re(s)>1$ we obtain that beautiful relation between $M, \zeta$
and $\theta$:
\begin{equation}
(M(p)\zeta)(s)\;=\;\int_{0}^{\infty}\theta(p)(x)x^{s-1}dx\;=\;M(\theta(p))
(s),
\end{equation}
since the above series is absolutely convergent (i.e. $p$ satisfies the
condition ($P_{2}$)) and we can interchange the order of summation and
integration.

Since $p \in \LP_{1}$, then $p \in Inv(\LR_{+})$ and satisfies $(P_{3})$
and ($P_{4}$), i.e. the following two series :
\begin{displaymath}
\sum_{n \in \LZ}p(x+ny)\;\;and\;\;\sum_{n \in \LZ}\hat{p}(x+ny)
\end{displaymath}
are {\bf absolutely} and {\bf almost uniformly convergent} for all $x, y
\in \LR_{+}$. In particular, such $p$ is $(\LR, \LZ)-admissible$, in the sense
of Weil's general assumption (see [We, VII.2]).
As the consequence the (PSF) holds. Using its and the initial condition
$(P_{0}) : p(0)=1$, changing variables, we can write
\begin{equation}
(M(p)\zeta)(s)\;=\;\int_{0}^{1}x^{s-2}\theta(p)(\frac{1}{x})dx
\;+\;\int_{1}^{\infty}x^{s-1}\theta(p)(x)dx\;=\;
\end{equation}
\begin{displaymath}
\;=\;\int_{1}^{\infty}[x^{-s}\theta(\hat{p})(x)\;+\;x^{s-1}\theta(p)(x)]dx\
;=:\;I(\theta(p))(s).
\end{displaymath}
The integral on the right-hand side of (3.39) converges uniformly for
$-\infty < a \le Re(s) \le b <+\infty$, since for $x \ge 1$, we have 
\begin{displaymath}
\mid x^{-s} \mid \le x^{-a} \;\;and\;\;\mid x^{s-1} \mid \le x^{b-1},
\end{displaymath}
while, for $l>1$ (obviously it suffices to take $l=2$ - like in ($P_{2}$
and $(P_{3})$)) holds
\begin{displaymath}
\theta(p)(x)\;\le\;\sum_{n=1}^{\infty}\mid p(nx) \mid\; \le
\;\frac{sup_{x \in \LR_{+}}\mid p(x) \mid x^{l}\zeta(l)}{x^{l}},
\end{displaymath}
and
\begin{displaymath}
\theta(\hat{p})(x)\;\le\;\frac{sup_{x \in \LR_{+}}\mid hat{p}(x) \mid
x^{l}\zeta(l)}{x^{l}}.
\end{displaymath}
Therefore, for each $p \in \LP_{1}$, the integral in (3.39) represents an
{\bf entire function} of $s$. Moreover - as it is well-known - $M(exp^{-1})(s)
= \Gamma(s)$ does not {\bf vanishes} anywhere (obviously $exp^{-1} \in
\LP_{1}$).

Thus the {\bf $\theta \Gamma$-quotient}
\begin{equation}
\zeta(s)\;:=\zeta(exp^{-1})(s)\;:=\;\frac{1}{s(s-1)\Gamma(s)}\;+\;
\frac{I(\theta(exp^{-1}))(s)}{\Gamma(s)}
\;,\;s \in \LC,
\end{equation}
gives the {\bf meromorphic continuation} of "local" zeta to the whole
complex plane.

If now $p \in \LP_{1}-\{exp^{-1}\}$ (let us observe that in this case
$M(p)$ {\bf can have zeros}!), then again acoording to (3.39) we have the
identity
\begin{equation}
(M(p)\zeta)(s)\;=\;\frac{1}{s(s-1)}\;+\;I(\theta(p))(s)\;\;,\;Re(s)>1.
\end{equation}
But now, the left-hand and right-hand sides of the above equation
(3.41) are the analytic functions in $D:= Re(s)>0 - \{1\}$. Thus, they must
be equal in $D$, according to the {\bf uniqueness} of the analytic
continuation of a holomorphic function in a domain (see e.g. [M,
XV.2]).

In the equations (3.36) - beside $\zeta, M(p), {\cal F}(p)$ and
$\theta(p), p \in \LP_{1}$ - appears the simple algebraic but very
important (deciding) meromorphic function $\frac{1}{s(s-1)}$ - from the
point of view of the Riemann hypothesis.

It is very convenient to look at this moment on $\frac{1}{(s-1)}$
(since $\zeta(s)$ has only the unique pole in $s=1$, so $s=0$ is only
an "apparent pole" of $\zeta$) as on the {\bf trivial zeta} of
$\LP^{1}(\LC)$, since its form is {\bf surprisely close} to the form of
the {\bf congruence Weil zeta} $Z(\LP^{1}(\tilde{\LF}_{q})/\LF_{q};s)$
of the {\bf projective line} $\LP^{1}(\tilde{\LF}_{q})$ over the
algebraic closure $\tilde{\LF}_{q}$ of the finite {\bf q-elements field}
$\LF_{q}$ (cf. e. g. [Kob, V.1]) :
$Z(\LP^{1}(\tilde{\LF}_{q})/\LF_{q};s) =\frac{1}{(1-s)(1-qs)}$. 
Since $q = p^{d}$ for some positive integer
$d$, where $p = char(\LF_{q})$ is the {\bf characteristic} of the field
$\LF_{q}$, and -as it is well-known, the characteristic of the field
$\LC : char(\LC) =q$ is zero : $q=0$. Thus formally : $\frac{1}{s-1} =:
Z(\LP^{1}(\tilde{\LF}_{0}/\LF_{0};s))$, i.e. we formally put
$\tilde{\LF}_{0} = \LF_{0} :=\LC$.

We define also the {\bf trivial Riemann zeta} $\zeta_{t}$ by the
trivial formula :
\begin{equation}
\zeta_{t}(s)\;:=\;Im(s)(2Re(s)\;-\;1)\;=\;\mid s(s-1)\mid^{2}
Im(\frac{1}{s(s-1)})\;,\;s\in \LC.
\end{equation}

\section{The Brownian motion and Fubini theorem as a sufficent and efficent
tools for a proof of the Riemann hypothesis.}
Let ${\cal M}$ be the {\bf category of all functional measure spaces}. Thus 
an object of ${\cal M}$ is a triple $(M, {\cal A}, \mu)$, where $M=F(\LR_{+})$ 
is a non-empty set of functions $f : \LR_{+}\longrightarrow \LC$, 
endowed with a $\sigma$-field ${\cal A}$ of subsets of $M$ and $\mu
: {\cal A}\longrightarrow \LR_{+}$  is a positive $\sigma$-additive
measure. Then, we can say on the {\bf measure characteristic numbers}
$m_{s} : {\cal M}\longrightarrow \LC\cup \{\infty\}$ as the family of 
{\bf functors on the category ${\cal M}$} indexed by $s \in \LC$ and defined
by the formula : let $(M, {\cal A}, \mu) \in {\cal M}$ be arbitrary. Then
\begin{equation}
m_{s}(M, {\cal A}, \mu)\;:=\;\int\int_{M \times
\LR_{+}}x^{s-1}f(x)d\mu(f)dx.
\end{equation}
Let us observe at once three elementary properties of the measure
characteristic numbers :

(1) $m_{s}$ is a {\bf measure invariant}, i.e. if two measure spaces
$(X, {\cal A}, \mu)$ and $(Y, {\cal B}, \nu)$ are {\bf isomorphic} in
the category ${\cal M}$ (we write then $(X, {\cal A}, \mu) \simeq_{M}
(Y, {\cal B}, \nu)$), i.e. there exists a set isomorphism $f :X
\longrightarrow Y$ with the property that ${\cal A} = f^{-1}({\cal B})$
and $\nu = f^{*}(\mu)$ ($\nu$ is the transport of $\mu$ by $f$ ), then 
\begin{equation}
m_{s}(X, {\cal A}, \mu)\;=\;m_{s}(Y, {\cal B}, \nu)\;\;,\,s \in \LC.
\end{equation}
\begin{equation}
(2)\;\;\mid m_{s}(\LP_{1}, {\cal P}, r_{m})\mid \;\le\;b_{s}(\LP_{1}, {\cal
P}, r_{m}),
\end{equation}
and finally $m_{s}$ are {\bf submiltiplicative}, i.e.
\begin{equation}
(3)\;\;m_{s}(X_{1}\times X_{2}, {\cal A}_{1}\otimes {\cal A}_{2},
\mu_{1}\times \mu_{2})\;\le m_{s}(X_{1}, {\cal A}_{1},
\mu_{1})m_{s}(X_{2}, {\cal A}_{2}, \mu_{2}).
\end{equation}
\begin{re}
Let us remark that the above defined measure and Brownian characteristic numbers
defined in the category ${\cal M}$ , are some analogs of the Betti numbers
and the Euler-Poincare characteristic in the category ${\cal T}$ of the
topological spaces. Like $m_{s}$ are measure invariants (see (1)) the
Euler-Poincare characteristic $\chi(X)$ of a topological space $X$ is
obviously the {\bf topological invariant}( even a {\bf homotopy invariant}).
Let us mention here its three fundamental properties:

(i) let $p : E \longrightarrow B$ be a locally trivial vector bundle with
the fibre $F$. Then (under some restrictions on the spaces $E, B, F$) their
characteristics are associated by the formula :$\chi(E) =\chi(B)\chi(F)$.

(ii) In particular, an Euler-Poincare characteristic of a direct product
of two topological spaces is equal their product, i.e. $\chi(X \times Y)=\chi
(X)\chi(Y)$.

(iii) With the help of the relation : $\chi(A \cup B)=\chi(A)+\chi(B)-\chi
(A \cap B)$, which is true for any cutting triple $(A \cup B, A, B)$, one can
calculate the Euler-Poincare characteristic of any 2-dimensional compact 
manifold. In particular $\chi(\cdot)$ is a finitely additive measure.
\end{re}

\begin{th}
{\bf The Riemann hypothesis is true}.
\end{th}
{\bf Proof}. According to Th.2, we have two real valued equations
\begin{equation}
Im[(M(p)\zeta)(s)]\;=\;\frac{\zeta_{t}(s)}{\mid s(s-1)
\mid^{2}}\;+\;\int_{1}^{\infty}[x^{Re(s)-1}\theta(p)(x)\;-\;x^{-Re(s)}
\theta(\hat{p})(x)]sin(Im(s)x)dx.
\end{equation}

We integrate (4.48) with respect to the {\bf Wiener-Riemann measure} $r_{m}$
and obtain :
\begin{equation}
Im[\int_{\LP_{1}}(M(p)\zeta)(s))dr_{m}(p)]\;=\;Z(\LP^{1}(\LC)/\LC;s)
\zeta_{t}(s)\;+\;
\end{equation}
\begin{displaymath}
\;+\;Im(\int_{\LP_{1}}dr_{m}(p)\int_{0}^{\infty}dx
\sum_{n=1}^{\infty}[x^{s-1}p(nx)\;+\;x^{-s}\hat{p}(nx)]).
\end{displaymath}
According to the Th.1 ($RH_{4}$) we also have
\begin{displaymath}
(\int_{\LP_{1}}dr_{m}(p)\int_{0}^{\infty}x^{s-1}p(x)dx)\zeta(s)=m_{s}(\LP_{1},
{\cal P}, r_{m})\zeta(s),
\end{displaymath}
where $m_{s}(\LP_{1}, {\cal P}, r_{m})$ is the measure characteristic
number of that measure space and obviously : $0 < Re(s) < 1/2$.

For the right-hand side of (4.50) we have an "easy" Fubini-Tonelli theorem.

Reelly, let us consider the triple iterated integral
\begin{equation}
I_{cxr}(\alpha)\;:=\;\sum_{n=1}^{\infty}\int_{1}^{\infty}x^{\alpha-1})dx
\int_{\LP_{1}}\mid p(nx) \mid dr_{m}(p).
\end{equation}
Then, making the substitution : $nx=t$ and using the property $(R_{4})$
of the Theorem 1 we get
\begin{equation}
I_{cxr}(\alpha)\;=\;(\sum_{n=1}^{\infty}\frac{1}{n^{2-\alpha}})\int_{n}
^{\infty}t^{\alpha-1}G(t)dt\int_{\LP_{1}} \mid p(t) \mid dr_{m}(p)\le
\zeta(2-\alpha)b_{\alpha}<+\infty,
\end{equation}
whereas $0<\alpha<1/2$.

Analogously, let us denote
\begin{equation}
I_{cxyr}(\beta):=\sum_{n=1}^{\infty}\int_{1}^{\infty}x^{-\beta}dx\int
_{\LR}dy \int_{\LP_{1}} \mid p(nxy) \mid
dr(p)=\sum_{n=1}^{\infty}\int_{1}^{\infty}x^{-\beta}dx \int_{\LR}G(nxy)
E \mid B^{0}_{\sqrt{nxy}}\mid dy \le
\end{equation}
\begin{displaymath}
\le \sum_{n=1}^{\infty}\int_{1}^{\infty}x^{-\beta}dx
\int_{\LR}G(nxy)(nxy)^{1/4}dy.
\end{displaymath}
Making the substitution : $n^{2}xy = z, dy = z/n^{2}x$ in the inner
integral on $\LR$, we get 
\begin{equation}
I_{cxyr}(\beta)\le
\sum_{n=1}^{\infty}\int_{1}^{\infty}\frac{dx}{n^{2}x^{1+\beta}}\int_{\LR}
G(z/n)(z/n)^{1/4} dz=\zeta(2+1/4)\int_{1}^{\infty}\frac{dx}{x^{1+\beta}}
\int_{\LR}G(z)z^{1/4}dz<+\infty ,
\end{equation}
for all $\beta >0$. 
According to the {\bf Tonelli theorem} all above triple iterated
integrals coincides with $I_{cxr}(\alpha)$ and they are {\bf finite}
for $\alpha \in (0, 1/2)$ whereas all changings of iteration in
$I_{cxyr}(\beta)$ coincides and are {\bf finite} for all $\beta>0$.

Thus finally, the averaging of the {\bf Muntz's relations} - given in
Th.2 - with respect to the {\bf Wiener-Riemann measure} $r_{m}$ leads
to the following functional equation : let $s = u + iv$. Then
\begin{equation}
Im[\zeta(s)\int_{0}^{\infty}x^{s-1}dx
\int_{\LP_{1}}p(x)dr_{m}(p)]\;=\frac{\int_{\LP_{1}}p(0)dr_{m}(p)\zeta_{t}
(s)}{\mid \zeta_{\LP^{1}}(s)\mid^{2}}\;+\;
\end{equation}
\begin{displaymath}
+\sum_{n=1}^{\infty}\int_{1}^{\infty}x^{u-1}sin(vlog x)dx \int_{\LP_{1}}p(nx)
dr_{m}(p)-\sum_{n=1}^{\infty}\int_{1}^{\infty}x^{-u}sin(vlog x)dx
\int_{\LR}cos(2\pi nxy)dy \int_{\LP_{1}}p^{+}(y)dr_{m}(p).
\end{displaymath}

According to the properties $(RH_{0})-(RH_{3}))$ of Th.1, we have
\begin{equation}
\int_{\LP_{1}}p(0)dr_{m}(p)=1\;,\;\int_{\LP_{1}}p(nx)dr_{m}(p)\;=\;0\;
if\;x \ge 1,
\end{equation}
and
\begin{equation}
\int_{\LP_{1}}p^{+}(y)dr_{m}(p) \;=\;m(\mid y \mid)\;if\;\mid y \mid
\le 1.
\end{equation}
Combining (4.55), (4.56) with (4.57) we can reduce the right-hand side of
(4.55) to the equation :
\begin{equation}
Im[m_{s}(\LP_{1}, {\cal P},r_{m})\zeta(s)]\;=
\;\frac{\zeta_{t}(s)}{\mid \zeta_{\LP^{1}}(s)\mid^{2}}-
\sum_{n=1}^{\infty}\int_{1}^{\infty}x^{-u}sin(v logx)\hat{m}(nx)dx.
\end{equation}
We have the following easy approximation for the triple integral which
appears in the right-hand side of (4.58) :
\begin{equation}
\mid \sum_{n=1}^{\infty}\int_{1}^{\infty}x^{-u}sin(v log x)dx
\int_{\LR}cos(2\pi nxy)m(y)dy \le
\sum_{n=1}^{\infty}\int_{1}^{\infty}x^{-u} \mid \hat{m}(nx) \mid dx.
\end{equation}
Since obviously $\hat{m} \in {\cal S}(\LR)$, then the right-hand side
of (4.59) is less then
\begin{equation}
(\sum_{n=1}^{\infty}\frac{1}{n^{2}})\int_{1}^{\infty}
\frac{x^{-u}dx}{x^{2}}max_{n \in \LN, x \ge 1}\mid
(nx)^{2}\hat{m}(nx)\mid \le
\zeta(2)\int_{1}^{\infty}\frac{dx}{x^{2+m}}x_{0}^{2} \mid m (x_{0}) \mid.
\end{equation}

Let us take any sequence $\{m_{k}\}$ of moment functions which
converges pointwisely to the {\bf characteristic function} $\chi_{0}(x)$
of $\{0\}$ in $\LR$ ( not to the Dirac delta distribution $\delta_{0}(x)$!).
Since $\mid m_{k}(x) \mid \le  \chi_{[0,1]}(x)$, using the {\bf Lebesgue's
dominated convergence theorem}, for $Re(s) \in (0,1/2)$ we finally get 
\begin{equation}
\lim_{k \rightarrow \infty}Im[m_{s}(\LP_{1}, {\cal P},
r_{m_{k}})\zeta(s)]\;=\;\zeta_{t}(s)\mid \zeta_{\LP^{1}(\LC)}(s) \mid^{2}.
\end{equation}
Since the non-trivial Riemann zeta $\zeta$ zeros lies symmetricaly with 
respect to the lines : (i) critical $Re(s)=1/2$ and (ii) $Im(s) =0$, then 
the equation (4.61) gives the most direct and an integral form proof of the 
{\bf Main Algebraic Conjecture} (MAC in short) :
\begin{displaymath}
(MAC)\;\;\zeta(\LC)\; \subset\; \zeta_{t}(\LC).
\end{displaymath}
Thus (RH) is proved.

\begin{re}
In [MA, Remark.20] we observed a strange violation of symetry with respect to 
the functions $Re(s)$ and $Im(s)$ in the Riemann hypothesis problem.
Roughly speaking it consists on the fact that the function $Im(s)$ is
much most important than $Re(s)$ (and in fact fundamental) for (RH).
That strategy also is explored in [MH], [MR], [MD] and [AM] (as well as
in the prepareted now paper [MRam]). In this paper we "broken" that
violation of the symmetry, i.e. in this paper $Im(s)$ and $Re(s)$ are
"equally important".
But in our opinion, the mentioned above violation of the symmetry was
one from reasons that (RH) was an open problem for 150 years. A trial
of a partial explanation of that phenomena, was presented in [MA]. Here
we give - maybe - a better explanation of that one : let us observe,
that instead of the case of $Re(s)$, the function $Im(s)$ seems to be
much more interested from the another point of view: it is a simple
example of the {\bf non-holomorphic modular form} of weight 2.

On the other hand - for some time it has been known - that (RH) is
strictly connected with the fundamental properties of the theta modular forms.
\end{re}

e-mails: madrecki@o2.pl, andrzej.madrecki@pwr.wroc.pl
\end{document}